\documentclass[a4paper,12pt]{amsart}

\usepackage{amsmath,graphics}
\usepackage{amssymb}
\usepackage{amsfonts}
\usepackage{latexsym}
\usepackage{mathrsfs}
\usepackage[dvips]{graphicx}
\usepackage{color}
\usepackage{enumerate}
\newtheorem{thm}{Theorem}[section]

\newtheorem{propo}[thm]{Proposition}
\newtheorem{lem}[thm]{Lemma}

\renewcommand{\Re}{{\rm Re}}
\renewcommand{\Im}{{\rm Im}}
\newcommand{\R}{\mathbb{R}}
\newcommand{\C}{\mathbb{C}}

\newcommand{\N}{\mathbb{N}}

\renewcommand{\H}{\mathbb{H}^2}

\newcommand{\lt}{{\mathscr L}}

\newcommand{\half}{{\textstyle{\frac{1}{2}}}}
\newcommand{\halfdelta}{{\textstyle{\frac{\delta}{2}}}}

\newcommand{\D}{{\mathcal D}}

\begin{document}

\bibliographystyle{plain}
\title[Density and localization of resonances]{Density and localization of resonances for convex co-compact hyperbolic surfaces}

\author[Fr\'ed\'eric Naud]{Fr\'ed\'eric Naud}
\address{
Laboratoire d'Analyse non-lin\'eaire et G\'eom\'etrie (EA 2151)\\
Universit\'e d'Avignon et des pays de Vaucluse, F-84018 Avignon, France.
}
\email{frederic.naud@univ-avignon.fr}

 \maketitle
 
\begin{abstract}
Let $X=\Gamma\backslash \H$ be a convex co-compact hyperbolic surface. We show that the density of resonances of the Laplacian $\Delta_X$ in strips 
$\{ \sigma \leq \Re(s)\leq \delta \}$  with $\vert \Im(s)\vert \leq T$ is less than $O(T^{1+\delta-\varepsilon(\sigma)})$ with $\varepsilon(\sigma)>0$ as long as $\sigma>\halfdelta$.
This improves the previous fractal Weyl upper bound of Zworski \cite{Zworski2} and is in agreement with the conjecture of \cite{JakobsonNaud2} on the essential spectral gap.
\end{abstract}
\begin{section}{Introduction and results}
\noindent
In this work, we will focus on the distribution of resonances of the Laplacian for a class of hyperbolic Riemann surface of infinite volume.
The spectral theory of these objects can be viewed as a relevant picture for more realistic physical models such as obstacle or potential scattering, but is also interesting in itself because
of its connection with counting problems and number theory (see for example the recent work of Bourgain-Gamburd-Sarnak \cite{BGS}, or Bourgain-Kontorovich \cite{BK1}).
Resonances replace the missing eigenvalues in non-compact situations and their precise distribution and localization in the complex plane is still a widely open subject. 
Let us be more specific. Let $\H$ denote the usual hyperbolic plane with its standard metric with curvature $-1$, and let $X=\Gamma \backslash \H$ be a convex co-compact hyperbolic surface (see next $\S$ for more details). Here $\Gamma$ is a discrete group of isometries of $\H$, convex co-compact (that is to say finitely generated, without non trivial parabolic elements or elliptic elements). In this paper we will assume that $\Gamma$ is {\it non-elementary} which is equivalent to say that $X$ is not a hyperbolic cylinder.
The limit set $\Lambda(\Gamma)$ is commonly defined as
$$\Lambda(\Gamma):=\overline{\Gamma.z}\cap \partial \H,$$
where $\Gamma.z$ denotes the orbit of $z\in \H$ under the action of $\Gamma$. This limit set actually does not depend on the choice of $z$ and
is a Cantor set in the above setting. We will denote by $\delta(\Gamma)$ the Hausdorff dimension of $\Lambda(\Gamma)$. Let $\Delta_X$ denote the hyperbolic
Laplacian on the surface $X$. The spectral theory on $L^2(X)$ has been described by 
Lax and Phillips \cite{LP2} :$[1/4, +\infty)$ is the continuous spectrum, has no embedded eigenvalues.
The rest of the spectrum is made of a (possibly empty) finite set of eigenvalues, starting at $\delta(1-\delta)$. The fact that the bottom of the spectrum is related to the dimension $\delta$ was first pointed out by Patterson \cite{Patterson1} for convex co-compact groups. This result was later extented for geometrically finite groups by Sullivan \cite{Sullivan,Sullivan2}.

\bigskip \noindent
By the preceding description of the spectrum, the resolvent 
$$R_X(s)=\left(\Delta_X-s(1-s)\right)^{-1}:L^2(X)\rightarrow L^2(X),$$
is therefore well defined and analytic on the half-plane $\{ \Re(s)>\half \}$ except at a possible finite set
of poles corresponding to the finite point spectrum. {\it Resonances} are then defined as poles of the meromorphic continuation of 
$$R_X(s):C_0^\infty(X)\rightarrow C^{\infty}(X)$$ 
to the whole complex plane. The set of poles is denoted by ${\mathcal R}_X$. This continuation is usually performed via the analytic Fredholm theorem
after the construction of an adequate parametrix. The first result of this kind in the more general setting
of asymptotically hyperbolic manifolds is due to Mazzeo and Melrose \cite{MazzMel}.
A more precise parametrix for surfaces was constructed by Guillop\'e and Zworski \cite{GuiZwor,GuiZwor1}. Note that by the above construction and choice of spectral parameter $s=\sigma+it$, the resonance set (including possible eigenvalues) ${\mathcal R}_X$ is included in the half-plane $\{ \Re(s)\leq \delta \}$. If $\delta>\half$, then except for a finite set of eigenvalues, $\mathcal{R}_X\subset \{ \Re(s)<\half \}$. 

\bigskip \noindent
One of the basic problems of the theory is to localize resonances with the largest real part, which are resonances who play a key role in various asymptotic
problems, including hyperbolic lattice point counting and wave asymptotics. Another central and related question is the existence of a {\it fractal Weyl law} when
counting resonances in strips. More precisely, in the papers by Zworski and Guillop\'e-Lin-Zworski, \cite{glz,Zworski2}  they prove the following. For all $\sigma \leq \delta$ set 
$$N(\sigma,T):=\#\{ z \in {\mathcal R}_X\ :\ \sigma\leq \Re(z) \leq \delta\ \mathrm{and}\ 0\leq \Im(z)  \leq T \},$$
then for all $\sigma$, one can find $C_\sigma$ such that for all $T\geq 0$, one has 
$$N(\sigma,T)\leq C_\sigma \left(1+T^{1+\delta} \right).$$
The first upper bound of this type involving a "fractal" dimension is due to Sj\"ostrand \cite{Sjostrand1} for potential scattering, see also \cite{SZduke}.
In \cite{glz}, the authors give some numerical evidence that the above estimate may be optimal provided $\sigma$ is small enough, but no rigorous results
so far have confirmed this conjecture, and the best existing lower bound is a sublinear omega estimate, see \cite{GuiZwor2}. On the other hand, it is conjectured 
in \cite{JakobsonNaud2} that for all $\epsilon>0$, there are only {\it finitely many resonances} in the half-plane 
$$\left \{ \Re(s)\geq \frac{\delta}{2}+\epsilon  \right \},$$
hence suggesting that one has to take $\sigma\leq \delta/2$ to observe the optimal rate of growth of $N(\sigma,T)$ as $T\rightarrow +\infty$.
In this paper, we show the following which is in agreement with the previous conjecture and provides a more precise picture in term of the distribution
of resonances.
\begin{thm} 
\label{mainresult}
Let $\Gamma$ be a non-elementary convex co-compact group as above. Then for all $\sigma \geq \halfdelta$, one can find $\tau(\sigma)\geq 0$ such that
as $T\rightarrow +\infty$,
$$N(\sigma,T)=O\left( T^{1+\tau(\sigma)}\right),$$
where $\tau(\halfdelta)=\delta$, $\tau(\sigma)<\delta$ for all $\halfdelta<\sigma$. Moreover, one can find $\sigma_0>\halfdelta$ such that
the map $\sigma \mapsto \tau(\sigma)$ is real-analytic, convex and strictly decreasing on $[\halfdelta,\sigma_0]$. 
In particular, its derivative at $\sigma=\halfdelta$ satisfies $\tau'(\halfdelta)<0$.
\end{thm}
\noindent
The main achievement of that paper is that $\tau(\sigma)<\delta$ for all $\sigma>\halfdelta$ and says that the $O(T^{1+\delta})$ upper bound is non-optimal
whenever $\sigma>\halfdelta$. From a physics point of view, this result is in agreement with the idea that the density of resonances has to "peak" in a strip close to "half of the classical escape rate" which is exactly $\halfdelta$ in our setting. See the paper \cite{Zworskip} for some numerical study that support this conjecture, especially figure 2. We also refer the reader to \cite{Nonnenmacher1} for a comprehensive survey on questions related to fractal Weyl laws and spectral gaps for various open (chaotic) quantum systems.
An "explicit" expression for $\tau(\sigma)$ is provided in the last section (see formula (\ref{taufunction})), involving the {\it topological pressure} of the Bowen-Series map.
The above result is of course relevant for $\sigma$ close enough to the critical value $\halfdelta$ since we know from a previous work of the author 
\cite{Naud2} that there always exists
a spectral gap i.e. one can find $\epsilon>0$ such that 
$${\mathcal R }_X \cap \{ \Re(s) \geq \delta-\epsilon \}=\{ \delta \}.$$ Unfortunately, this $\epsilon$ is hardly explicit. Of course if $\delta>\half$, then this spectral gap is already given by Lax-Phillips theory, but the above theorem is still
meaningful since $\halfdelta<\half$.  

This result can also be compared with the spectral deviations obtained by N. Anantharaman \cite{NaliniWave} 
for the resonances  of the damped wave equation on negativey curved
manifolds, except that in our case no Weyl law is known rigorously. The techniques we rely on to prove Theorem \ref{mainresult} are therefore very different although
they also involve some ergodic theory and thermodynamical formalism. 
We would like to point out that Theorem \ref{mainresult} has a natural interpretation in terms of Selberg zeta function. Let us denote by $\mathcal P$ the
set of primitive closed geodesics on $X=\Gamma \backslash \H$ and given ${\mathcal C } \in {\mathcal P}$ we denote by $l(\mathcal{C})$ its length.
Selberg zeta function $Z_\Gamma(s)$ is usually defined through the infinite product ($\Re(s)$ is taken large enough)
$$Z_\Gamma(s):=\prod_{(\mathcal{C},k)\in \mathcal{P}\times \N}\left(1-e^{-(s+k)l(\mathcal{C})} \right).$$
It is known since the work of Patterson-Perry \cite{PatPer} that this function extends analytically to $\C$ and the non-trivial (non topological) zeros of
$Z_\Gamma(s)$ are the resonances. Therefore Theorem \ref{mainresult} can be read as a statement on zeros of Selberg's zeta function, and this is actually how we prove it. Let us make a few comments on the organization of the paper. In the next section we recall the transfer operator approach for Selberg's zeta function.
In section $\S 3$ we prove the necessary a priori bounds to control the growth of a (modified ) Selberg zeta function in strips, generalizing slightly  the method of
\cite{glz}. In $\S4$ we use a classical lemma of Littlewood to relate the counting function for resonances to a "mean square estimate". The idea of $\S6$, 
which is the core of that paper, is to exhibit some cancellations in the quadratic sums to beat the pointwise
bound of \cite{glz}. To achieve this goal some lower bounds on the derivatives of "off-diagonal phases" are required and proved in $\S 5$. The technique we use in $\S 6$ bears some similarity with the work of Dolgopyat \cite{Dolgopyat, Naud2} and can be viewed as an "averaged" Dolgopyat estimate.  
We would like to mention that the ideas presented here should extend without major difficulties to higher dimensional Schottky groups, actually only $\S 5$ needs to be significantly modified. 
\end{section}
\begin{section}{Transfer operator and Selberg's zeta function}
We use the notations of $\S 1$. Let $\H$ denote the Poincar\'e upper half-plane $\H=\{ x+iy\in \C\ :\ y>0\}$ endowed with its standard metric of constant $-1$ curvature
$$ds^2=\frac{dx^2+dy^2}{y}.$$ The group of (positive) isometries of $\H$ is naturally isomorphic to $\mathrm{PSL}_2(\R)$ through the action of 
$2\times 2$ matrices viewed as M\"obius transforms
$$z\mapsto \frac{az+b}{cz+d},\ ad-bc=1.$$ 
A Fuchsian Schottky group is a discrete subgroup of $\mathrm{PSL}_2(\R)$ built as follows. Let $\D_1,\ldots, \D_p,\D_{p+1},\ldots, \D_{2p}$
be $2p$ Euclidean {\it open} discs in $\C$ orthogonal to the line $\R\simeq \partial \H$. We assume that for all $i\neq j$, $\D_i \cap \D_j=\emptyset$. 
Let $\gamma_1,\ldots,\gamma_p \in \mathrm{PSL}_2(\R)$ be $p$ isometries such that for all $i=1,\ldots,p$, we have
$$\gamma_i(\D_i)=\widehat{\C}\setminus \overline{\D_{p+i}},$$
where $\widehat{\C}:=\C\cup \{ \infty \}$ stands for the Riemann sphere. 



\bigskip \noindent
The discrete group $\Gamma$ generated by $\gamma_1,\ldots,\gamma_p$ and their inverses
is called a {\it classical Schottky group}. If $p>1$ then $\Gamma$ is said to be non elementary.
It is always a free, geometrically finite, discrete group and if in addition we require that for all $i\neq j$,
$\overline{\D_i}\cap \overline{\D_j}=\emptyset$, then $\Gamma$ is a convex co-compact group i.e. the quotient Riemann surface
$$X=\Gamma \backslash \H$$
is an {\it infinite volume} geometrically finite hyperbolic surface with no cusps. The converse is true, up to an isometry, all convex co-compact hyperbolic surfaces
can be uniformized by a group as above, see \cite{button}.

Let $\Gamma\subset \mathrm{PSL}_2(\R)$ be a Fuchsian Schottky group as defined earlier:
\[\Gamma=\langle \gamma_1,\ldots,\gamma_p;\gamma_1^{-1},\ldots,\gamma_p^{-1} \rangle, \]
where $\gamma_i(\D_i)=\widehat{\C}\setminus \overline{\D_{p+i}}$. We also set for $i=1,\ldots,p$, $\gamma_{p+i}:=\gamma_i^{-1}$. 
For all $j=1,\ldots,2p$, let $H^2(\mathcal{D}_j)$ denote the Bergman space of holomorphic functions defined by
\[ H^2(\mathcal{D}_j):=\left \{ f:\mathcal{D}_j\rightarrow \C\ :\ f\ \mathrm{holomorphic\ and}\ \int_{\mathcal{D}_j} \vert f \vert^2 dm<+\infty \right \}, \]
here $m$ stands for the usual Lebesgue measure.
Each function space $H^2(\mathcal{D}_j)$ is a Hilbert space when endowed with the obvious norm. We set 
\[H^2:=\bigoplus_{j=1}^{2p}  H^2(\mathcal{D}_j). \]
The Ruelle transfer operator is a bounded linear operator $\lt_s:H^2\rightarrow H^2$ defined by ($z \in \mathcal{D}_i$, $s\in \C$ )
\[ (\lt_s(f))_i(z):=\sum_{j\neq i} (\gamma_j'(z))^s f_{j+p}(\gamma_j(z)), \]
with the notation $f=(f_1,\ldots,f_{2p})$, and $j+p$ is understood mod $2p$. 

\noindent
We have to say a few words about the complex powers here: we have $\gamma'_j(\D_i)\subset \C\setminus \R^-$ for $i\neq j$ so 
$\gamma'_j(z)^s$ is understood as
$$\gamma'_j(z)^s:=e^{s\mathbb{L}(\gamma'_j(z))},$$
where $\mathbb{L}$ is a complex logarithm on $\C\setminus \R^-$ which coincides on $\R^+\setminus\{0\}$ with the usual logarithm. For example,
one can take
$$\mathbb{L}(z)=\int_1^z \frac{d\zeta}{\zeta}.$$
This operator $\lt_s$ acts as a compact, trace class operator on $H^2$ and we refer the reader to \cite{Borthwick, glz} for a proof. The Fredholm determinant 
associated to this family is the Selberg zeta function defined above: for all $s\in \C$, we have 
\[ Z_\Gamma(s)=\mathrm{det}(I-\lt_s).\]
This remarkable formula that dates back to the work of Pollicott \cite{Pollvap} in the co-compact case, is the starting point of our analysis. 
Unfortunately, the space $H^2$ as defined above is not enough for our analysis and we will need to use as in \cite{glz} a family of Hilbert spaces $H^2(h)$ 
depending on a small scale parameter $0<h\leq h_0$. We recall the construction taken from \cite{glz}, see also \cite{Borthwick} chapter 15. Let $0<h$ and
set
$$\Lambda(h):=\Lambda+(-h,+h),$$
then for all $h$ small enough, $\Lambda(h)$ is a bounded subset of $\R$ whose connected components have length at most $Ch$ where $C>0$ is independent of $h$, see \cite{Borthwick} Lemma 15.12. Let $I_\ell(h)$ denote these connected components, with $\ell=1,\ldots,N(h)$. The existence of a finite Patterson-Sullivan measure 
$\mu$ supported by $\Lambda(\Gamma)$ plus Sullivan Shadow Lemma (see \cite{Borthwick}, chapter 14) show that
$$ A^{-1}h^\delta N(h)\leq \sum_\ell \mu(I_\ell(h))=\mu(\Lambda(\Gamma)),$$
for some uniform $A>0$, hence the number $N(h)$ of connected components is $O\left(h^{-\delta} \right)$. Given $1\leq \ell \leq N(h)$, let $\D_\ell(h)$ be the unique euclidean open disc
in $\C$ orthogonal to $\R$ such that $$\D_\ell(h)\cap \R=I_\ell(h).$$ 
Now set 
$$H^2(h):= \bigoplus_{\ell=1}^{N(h)}H^2(\D_\ell(h)).$$ 
We will see in the next section that for $n\geq n_0$ (with $n_0$ independent of $h$) that the operators $\lt_s^n$ act as compact trace class operators
on $H^2(h)$. Moreover, the Fredholm determinant 
$$Z_\Gamma^{(n)}(s):=\det(I-\lt_s^n)$$
is a multiple of the Selberg zeta function and we will count resonances using this determinant instead of the original
zeta function. The goal of the next section is to provide a proof of the following fact, which is a slight modification of the argument of 
Guillop\'e-Lin-Zworski in \cite{glz}.
\begin{propo}
\label{est1}
For all $\sigma_0 \leq \delta$, there exists $C>0$ such that for all $\sigma_0\leq \Re(s) \leq \delta$ and $n\geq n_0$, we have for $\vert \Im(s) \vert$ large,
$$\log\vert Z_\Gamma^{(n)}(s) \vert \leq C \vert \Im(s) \vert^{\delta}e^{nP(\sigma_0)},$$
where $P(\sigma)$ is the topological pressure at $\sigma$.
\end{propo}
\noindent
We refer the reader to the next $\S$ for a definition of topological pressure.
\end{section}
\begin{section}{Basic pointwise estimates }
The goal of this section is to prove the previous Proposition.
We first introduce some notations. We recall that $\gamma_1,\ldots,\gamma_p$ are generators of the Schottky group $\Gamma$.
Considering a finite sequence $\alpha$ with
\[\alpha=(\alpha_1,\ldots,\alpha_n)\in \{1,\ldots, 2p\}^n,\]
we set 
\[ \gamma_\alpha:=\gamma_{\alpha_1}\circ \ldots \circ \gamma_{\alpha_n}. \]
We then denote by $\mathscr{W}_n$ the set of admissible sequences of length $n$ by
\[ \mathscr{W}_n:=\{ \alpha \in \{1,\ldots, 2p\}^n\ :\ 
\forall\ i=1,\ldots,n-1,\ \alpha_{i+1}\neq \alpha_i +p\ \mathrm{mod}\ 2p \}.\]
We point out that if $\alpha \in \mathscr{W}_n$, then $\gamma_\alpha$ is a {\it reduced word}
in the free group $\Gamma$. For all $j=1,\ldots, 2p$, we define $\mathscr{W}_n^j$ by
\[ \mathscr{W}_n^j:=\{ \alpha \in \mathscr{W}_n\ :\ \alpha_n\neq j \}. \] 
If $\alpha \in \mathscr{W}_n^j$, then $\gamma_\alpha$ maps $\overline{\D_j}$ into $\D_{\alpha_1+p}$.
Given the above notations and $f\in H^2$, we have for all $z\in \D_j$ and $n\in \N$,
$$\lt_s^n(f)(z)=\sum_{\alpha \in \mathscr{W}_n^j} (\gamma_\alpha'(z))^s f(\gamma_\alpha (z)).$$

We will need throughout the paper some distortion estimates
for these maps $\gamma_\alpha$. More precisely we have for all $j=1,\ldots,2p$,
\begin{itemize}
 \item {\it(Uniform hyperbolicity)}. One can find $C>0$ and $0<\overline{\theta}<\theta<1$ such that for all $n,j$ and $\alpha \in \mathscr{W}_n^j$,
 \[ C^{-1}\overline{\theta}^n\leq \sup_{\D_j}\vert \gamma'_\alpha \vert \leq C \theta^n. \]
 \item {\it (Bounded distortion).} There exists $M_1>0$ such that for all $n,j$ and all $\alpha \in \mathscr{W}_n^j$,
 \[ \sup_{\D_j} \left \vert \frac{\gamma''_\alpha}{\gamma'_\alpha}  \right \vert \leq M_1.\]
\end{itemize}
The second estimate, called bounded distortion, will be used constantly throughout the paper. In particular it implies
that for all $z_1,z_2 \in \D_j$, for all $\alpha \in \mathscr{W}_n^j$, we have
$$e^{-\vert z_1-z_2\vert M_1}\leq \frac{\vert \gamma'_\alpha(z_1) \vert}{\vert \gamma'_\alpha(z_2)\vert}\leq  e^{\vert z_1-z_2\vert M_1}.$$
These estimates are rather standard facts in the classical ergodic theory of uniformly expanding Markov maps. For a proof of the first estimate,
we refer the reader to \cite{Borthwick} for example. Bounded distortion follows from hyperbolicity and elementary computations.  
Another critical tool in our analysis is the {\it Topological pressure} and {\it Bowen's formula}. Let $I_j:=\D_j\cap \R$. The Bowen-Series map 
$T:\cup_{i=1}^{2p} I_i\rightarrow \R\cup\{\infty \}$ is defined by $T(x)=\gamma_i(x)$ if $x\in I_i$. The non-wandering set of this map is exactly the limit set 
$\Lambda(\Gamma)$ of the group:
\[ \Lambda(\Gamma)=\bigcap_{n=1}^{+\infty} T^{-n}(\cup_{i=1}^{2p} I_i).\]
The limit set is $T$-invariant and given a continuous map $\varphi:\Lambda(\Gamma)\rightarrow \R$, the topological pressure $P(\varphi)$ can be defined through the variational formula:
\[ P(\varphi)=\sup_{\mu}\left ( h_\mu(T)+\int_{\Lambda} \varphi d\mu \right),\]
where the supremum is taken over all $T$-invariant probability measures on $\Lambda$, and $h_\mu(T)$ stands for the measure-theoretic entropy.
A celebrated result of Bowen \cite{Bowen1} says that if one considers the map $\sigma\mapsto P(-\sigma \log \vert T' \vert)$ then this map is convex, strictly decreasing and vanishes
exactly at $\sigma=\delta(\Gamma)$, the Hausdorff dimension of the limit set. An alternative way to compute the topological pressure is to look at weighted sums
on periodic orbits i.e. we have 
\begin{equation}
\label{pressure1}
e^{P(\varphi)}=\lim_{n\rightarrow +\infty} \left( \sum_{T^n x=x} e^{\varphi^{(n)}(x)}\right)^{1/n},
\end{equation}
with the notation $\varphi^{(n)}(x)=\varphi(x)+\varphi(Tx)+\ldots +\varphi(T^{n-1}x).$
We will seldom use the two previous formulas in our analysis but will rather use the following upper bound. For simplicity, we will use the notation
$P(\sigma)$ in place of $P(-\sigma \log \vert T' \vert)$.
\begin{lem} For all $\sigma_0<M$ in $\R$, one can find $C_0>0$ such that for all $n\geq 1$ and $M\geq \sigma\geq \sigma_0$, we have
\begin{equation}
\label{pressure2}
 \sum_{j=1}^{2p}\left (\sum_{\alpha \in \mathscr{W}_n^j} \sup_{I_j}  \vert \gamma'_\alpha \vert ^\sigma\right)\leq C_0 e^{nP(\sigma_0)}.
\end{equation}
\end{lem}
\noindent {\it Proof.} For all $\sigma\geq \sigma_0$ and $x \in I_j$, we have by bounded distortion property
$$\sum_{\alpha \in \mathscr{W}_n^j} \sup_{I_j}  \vert \gamma'_\alpha \vert ^\sigma \leq 
C^M \sum_{\alpha \in \mathscr{W}_n^j}  (\gamma'_\alpha(x) )^{\sigma_0}=C^M \lt_{\sigma_0}^n(\mathbf{1})(x).$$
The Ruelle-Perron-Frobenius theorem \cite{Baladi} applied to $\lt_{\sigma_0}:C^1(\overline{I})\rightarrow C^1(\overline{I})$ ($I=\cup_j I_j$) says that
this operator is quasi-compact, has spectral radius $e^{P(\sigma_0)}$ and $e^{P(\sigma_0)}$ is the only eigenvalue on the
circle 
$$\{ \vert z \vert= e^{P(\sigma_0)} \}.$$
Moreover, $e^{P(\sigma_0)}$ is a simple eigenvalue whose  spectral projection is given by 
$$f\mapsto \mathcal{P}(f):=\varphi_0 \int_I f d\mu_0,$$
where $\mu_0$ is a probability measure and $\varphi_0$ is a positive $C^1$ density on $I$. Using Holomorphic functional calculus,
we can therefore decompose
$$\lt_{\sigma_0}^n=e^{nP(\sigma_0)}\mathcal{P}+\mathcal{N}^n, $$
where $\mathcal{N}$ has spectral radius at most $\theta_0 e^{P(\sigma_0)}$ for some $\theta_0<1$. It is now clear that
we have
$$\sum_{\alpha \in \mathscr{W}_n^j} \sup_{I_j}  \vert \gamma'_\alpha \vert ^\sigma \leq C_0 \sup \varphi_0 e^{nP(\sigma_0)}+C_1\widetilde{\theta_0}^n e^{nP(\sigma_0)},$$
for some constants $C_0,C_1$ and $\theta_0<\widetilde{\theta_0}<1$ and the proof is done. $\square$

\bigskip
\noindent We can now go back to the action of
$$\lt_s^n:H^2(h)\rightarrow H^2(h).$$
Set $\mathcal{E}_j(h):=\{ 1\leq \ell \leq N(h)\ :\ \D_\ell(h)\subset \D_j \}$.  Given $\alpha \in \mathscr{W}_n^j$ and $\ell \in \mathcal{E}_j(h)$, we have 
provided $n$ is large enough,
$$\gamma_\alpha(\D_\ell(h))\subset \D_{\ell'}(h)$$
for some $\ell'\in \{1,\ldots,N(h) \}$. Indeed, since the limit set $\Lambda$ is $\Gamma$-invariant, $\gamma_\alpha$ maps $\D_j(h)$ into a disc of size
at most $C h \theta^n$ which contains an element of $\Lambda$. Taking $n$ large enough (uniformly on $h$), this disc has to be inside $\cup_\ell \D_\ell(h)$.
Therefore $\lt_s^n$ leaves $H^2(h)$ invariant as long as $n$ is large. We will actually need something more quantitative which is proved in \cite{glz,Borthwick}:
\begin{lem}
 \label{separation}
 There exists $n_0$, $\kappa>0$ such that for all $n\geq n_0$, $\alpha\in \mathscr{W}_n^j$, we have for all 
 $\ell \in \mathscr{E}_j(h)$, $\gamma_\alpha(\D_\ell(h))\subset \D_{\ell'}(h)$ with
 $$ \mathrm{dist}(\gamma_\alpha(\D_\ell(h)),\partial \D_{\ell'}(h) )\geq \kappa h.$$
\end{lem}
\noindent
To estimate the pointwise growth of $Z_\Gamma^{(n)}(s)=\det(I-\lt_s^n)$ we will use some singular values estimates on the space $H^2(h)$ with $h=\vert \Im(s) \vert^{-1}$. Notice that this zeta function does not depend on $h$. Indeed, provided $\Re(s)>\delta$, we have the formula
\begin{equation}
\label{traceformula}
\det(I-\lt_s^n)=\exp\left(-\sum_{k=1}^\infty \frac{1}{k} \mathrm{Tr}( \lt_s^{nk})  \right),
\end{equation}
and the trace can be computed (using a fixed point analysis, see \cite{Borthwick}):
$$\mathrm{Tr}(\lt_s^p)=\sum_{T^px=x} \frac{(T^p)'(x)^{-s}}{1-((T^P)'(x))^{-1}},$$
the above sums over periodic orbits of the Bowen-Series are clearly independent of the choice of function space, hence of $h$ which will be adjusted to optimize our estimates.
First we need to recall some facts on singular values of compact operators and our basic reference is the book of Simon \cite{Simon}. If $\mathcal{H}$ is a complex Hilbert space and $T:\mathcal{H}\rightarrow \mathcal{H}$
is a compact operator, the singular values $\mu_k(T)$ are the eigenvalues (ordered decreasingly) of the self-adjoint positive operator $\sqrt{T^* T}$.
If we have
$$\sum_{k=1}^\infty \mu_k(T)<+\infty, $$
then $T$ is said to be trace class. If $T$ is trace class, then the eigenvalue sequence
$$\vert \lambda_1(T) \vert \geq \vert \lambda_2(T) \vert \geq \ldots \geq \vert \lambda_k(T)\vert$$
is summable and the trace $\mathrm{Tr}(T)$ is defined by
$$\mathrm{Tr}(T)=\sum_{k=1}^\infty \lambda_k(T),$$
while the determinant $\det(I+T)$ is given by
$$\det(I+T):=\prod_{k=1}^\infty \left( 1+\lambda_k(T)\right).$$
Weyl's inequalities show that we have
$$  \prod_{k=1}^\infty \left( 1+\vert \lambda_k(T)\vert \right)  \leq \prod_{k=1}^\infty \left( 1+\mu_k(T)\right).$$
In this paper we will use the following remark:
\begin{lem}
\label{det} 
Let $T:\mathcal{H}\rightarrow \mathcal{H}$ be a trace class operator and $(e_\ell)_{\ell \in \N}$ is a Hilbert basis of $\mathcal{H}$, then we have
$$\log\vert \det (I+T)\vert\leq \sum_{\ell=0}^\infty \Vert T(e_\ell) \Vert.$$
\end{lem}
\noindent {\it Proof}. By Weyl's inequality, we write
$$\log\vert \det (I+T)\vert\leq\sum_{k=0}^\infty \log(1+\mu_k(T))\leq \sum_{k=0}^\infty \mu_k(T)=\mathrm{Tr}(\sqrt{T^*T}).$$
But the classical Lidskii theorem says that
$$\mathrm{Tr}(\sqrt{T^*T})=\sum_{\ell=0}^\infty \langle \sqrt{T^*T} e_\ell,e_\ell \rangle \leq 
\sum_{\ell=0}^\infty \Vert \sqrt{T^*T} e_\ell \Vert,$$
by Schwartz inequality, but we have $\Vert \sqrt{T^*T} e_\ell \Vert=\Vert T e_\ell \Vert$. $\square$

\bigskip \noindent
We can now give a proof of Proposition \ref{est1}. First we need an explicit basis of $H^2(h)$. Setting for $\ell \in \{1,\ldots, N(h)\}$,
$$D_\ell(h):=D(c_\ell,r_\ell)$$
with $c_\ell\in \R$ and $h\leq r_\ell<Ch$, we denote by $\mathbf{e}_k^\ell$ the function in $H^2(h)$ defined for $z\in \D_j(h)$ by
$$\mathbf{e}_k^\ell(z)=\left \{ \begin{array}{c}
0\ \mathrm{if}\ j\neq \ell \\
\sqrt{\frac{k+1}{\pi}}\frac{1}{r_j} \left (\frac{z-c_j}{r_j} \right)^k\ \mathrm{if}\ j=\ell.
\end{array}\right. $$
It is straightworfard to check that $(\mathbf{e}_k^\ell)_{k \in \N,1\leq \ell \leq N(h)}$ is a Hilbert basis of $H^2(h)$. We now use Lemma \ref{det} and write
$$\log\vert Z_\Gamma^{(n)}(s)\vert\leq \sum_{k=0}^{\infty} \sum_{\ell=1}^{N(h)}\Vert \lt_s^{n}( \mathbf{e}_k^\ell )\Vert_{H^2(h)}.$$
Using Schwartz inequality, we get ($N(h)=O(h^{-\delta})$)
\begin{equation}
\label{ineq2}
\log\vert Z_\Gamma^{(n)}(s)\vert\leq Ch^{-\halfdelta} \sum_{k=0}^{\infty} \left (\sum_{\ell=1}^{N(h)}\Vert \lt_s^{n}( \mathbf{e}_k^\ell )\Vert_{H^2(h)}^2 \right)^{1/2}.
\end{equation}
On the other hand we have 
$$\sum_{\ell=1}^{N(h)}\Vert \lt_s^{n}( \mathbf{e}_k^\ell )\Vert_{H^2(h)}^2=
\sum_{j=1}^{2p}\sum_{\ell' \in \mathscr{E}_j(h)} \int_{\D_{\ell'}(h)}\sum_{\ell=1}^{N(h)}\vert \lt_s^n(\mathbf{e}_k^{\ell}) \vert^2 dm,$$
where $m$ is the usual Lebesgue measure on $\C$. Given $\ell' \in \mathscr{E}_j(h)$,  we can write
\begin{equation}
\label{ineq1}
 \int_{\D_{\ell'}(h)}\sum_{\ell=1}^{N(h)}\vert \lt_s^n(\mathbf{e}_k^{\ell}) \vert^2 dm
\leq  \sum_{\alpha,\beta \in \mathscr{W}_n^j} \int_{\D_{\ell'}(h)} \vert (\gamma'_\alpha)^s \vert \vert (\gamma'_\beta )^s \vert F_{\alpha,\beta}^{(k)} dm,
\end{equation}
where 
$$F_{\alpha,\beta}^{(k)}(z):=\sum_{\ell=1}^{N(h)} \vert \mathbf{e}_k^\ell \circ \gamma_\alpha(z) \vert  \vert \mathbf{e}_k^\ell \circ \gamma_\beta (z) \vert.$$
We now need to prove the following remark.
\begin{lem}
\label{baseestimate} Setting $\Omega_j(h):=\cup_{\ell \in \mathscr{E}_j(h)} \D_\ell(h)$, we have for all $\alpha,\beta \in \mathscr{W}_n^j$,
$$\sup_{\Omega_j(h)}\vert F_{\alpha,\beta}^{(k)}\vert\leq Ch^{-2}\rho^k,$$
with $C$ and $0<\rho<1$ uniform.
\end{lem}
\noindent {\it Proof}. First remark that given $z\in \Omega_j(h)$, we have either $F_{\alpha,\beta}^{(k)}(z)=0$ or
$$F_{\alpha,\beta}^{(k)}(z)=\vert \mathbf{e}_k^{\ell_0} \circ \gamma_\alpha(z) \vert  \vert \mathbf{e}_k^{\ell_0} \circ \gamma_\beta (z) \vert,$$
for some $\ell_0 \in \{1,\ldots,N(h)\}$. Then combine Lemma \ref{separation} with the explicit formula for $\mathbf{e}_k^{\ell_0}$ to obtain the result. $\square$

\bigskip \noindent Going back to estimate (\ref{ineq1}), and using the fact (use bounded distortion) that we set $\vert \Im(s) \vert=h^{-1}$
$$ \sup_{z\in \D_{\ell'}(h)}\vert( \gamma_\alpha')^s\vert \leq  C\Vert \gamma_\alpha '\Vert_{\infty,j}^{\Re(s)}, $$
we have reached
$$\sum_{\ell=1}^{N(h)}\Vert \lt_s^{n}( \mathbf{e}_k^\ell )\Vert_{H^2(h)}^2\leq C h^{-\delta} \rho^k
\left (\sum_{j=1}^{2p}\sum_{\alpha \in \mathscr{W}_n^j} \sup_{I_j}  \vert \gamma'_\alpha \vert ^\sigma\right)^2$$
$$\leq Ch^{-\delta} \rho^k e^{2nP(\sigma)}.$$
The proof is done by inserting the above estimate in formula (\ref{ineq2}) and summing over $k$. $\square$

 \end{section}
 \begin{section}{Applying Littlewood's Lemma}
 We now show how to reduce the proof of Theorem \ref{mainresult} to a mean square estimate on transfer operators. To this end we apply a
 result of Littlewood taken from the classics. More precisely, we prove the following. Define $M(\sigma,T)$ by
 $$M(\sigma,T):=\#\{ s\in \mathcal{R}_X\ :\ \sigma\leq \Re(s)\leq \delta\ \mathrm{and}\ T/2\leq \Im(s)\leq T\}.$$
 \begin{propo}
 \label{Little}
 Let $\sigma_0<\sigma<\delta$, then for all $T$ large and $n(T)=[\nu \log T ]$ with $\nu>0$ small enough, we have
 $$M(\sigma,T)\leq C_1\int_{T/2}^T\log\vert Z_\Gamma^{(n(T))}(\sigma_0+it)\vert dt+C_2 T .$$
 \end{propo}
 \noindent {\it Proof.} We start by recalling a version of Littlewood's Lemma which suits our needs. Let $\sigma_0<1$ and $T>0$. Let $f$ be a 
 function which is holomorphic on a neighborhood of the rectangle 
 $$R_{T,\sigma_0}:=[\sigma_0,1]+i[T/2,T],$$
 and assume that $f$ does not vanish on the segment $1+i[T/2,T]$. Denote by $\mathcal{Z}_{T,\sigma_0}$ the set of zeros of $f$ on $R_{T,\sigma_0}$.
 Then we have the formula
 $$2\pi \sum_{z\in \mathcal{Z}_{T,\sigma_0}} \left(\Re(z)-\sigma_0\right)=\int_{T/2}^T\log \vert f(\sigma_0+it)\vert dt-
 \int_{T/2}^T\log \vert f(1+it)\vert dt $$
 $$+\int_{\sigma_0}^1 \mathrm{Arg}(f(\sigma+iT))d\sigma- \int_{\sigma_0}^1 \mathrm{Arg}(f(\sigma+iT/2))d\sigma.$$
 The function $\mathrm{Arg}(f)$ is the imaginary part of a determination of a complex logarithm $\log f $ which is defined by taking upper limits of a holomorphic
 logarithm on a suitable simply connected domain (see Titchmarsh \cite{Tit2}, section 9.9, for more details). The resulting function $\mathrm{Arg}(f)$ is well defined on  $R_{T,\sigma_0}\setminus \mathcal{Z}$ and is discontinuous on a finite union of segments. 
 
 \bigskip \noindent Fix $\sigma_0<\sigma<\delta<1$. Then since $\mathcal{R}_X$ is a subset of the set of zeros of $Z_\Gamma^{(n)}(s)$,
 applying the above formula to $Z_\Gamma^{(n)}(s)$ we get 
 $$2\pi (\sigma-\sigma_0)M(\sigma,T)\leq \int_{T/2}^T\log\vert Z_\Gamma^{(n(T))}(\sigma_0+it)\vert dt$$
$$+ O\left( \int_{T/2}^T\log \vert Z_\Gamma^{(n(T))}(1+it)\vert dt \right)+O\left(\sup_{\sigma_0\leq \sigma\leq 1}\vert \mathrm{Arg}(Z_\Gamma^{(n(T))}(\sigma+iT))\vert \right)$$
$$+ O\left(\sup_{\sigma_0\leq \sigma\leq 1}\vert \mathrm{Arg}(Z_\Gamma^{(n(T))}(\sigma+iT/2))\vert \right).$$
Remark that combining (\ref{pressure1}) and the trace formula (\ref{traceformula}) shows ($P(1)<0$) that $Z_\Gamma^{n(T)}(s)$ is uniformly bounded
on $\{ \Re(s)=1 \}$, hence
$$\int_{T/2}^T\log \vert Z_\Gamma^{(n(T))}(1+it)\vert dt=O(T).$$ 
To control $\mathrm{Arg}(Z_\Gamma^{(n)}(s))$, we will use another classical result of Titchmarsh essentially based on Jensen's formula, see \cite{Tit2}, 9.4.
\begin{lem}
 Fix $2>\sigma_0>0$ and let $f$ be a holomorphic function on the half-plane $\{\Re(s)>0\}$. Suppose that
 $\vert \Re(f(2+it))\vert\geq m>0$ for all $t\in \R$ and assume that $f(\overline{s})=\overline{f(s)}$ for all $s$.
 Suppose in addition that $\vert f(\sigma+it)\vert\leq A_{\sigma,t}$, then if $T$ is not the imaginary part of a zero of $f(s)$, for all $\sigma\geq \sigma_0$,
 $$\vert \mathrm{Arg}(f(\sigma+iT))\vert \leq C_{\sigma_0}( \log A_{\sigma_0,T+2}-\log m)+3\pi/2.$$
\end{lem}
\noindent
To be able to apply the above Lemma, we need to prove the following.
\begin{lem}
There exists $m>0$ such that for all $n$ large enough, we have for all $t$,
$$\vert \Re(Z_\Gamma^{(n)}(2+it))\vert\geq m.$$
\end{lem}
\noindent {\it Proof}. First remark that for $\Re(s)>\delta$, 
$$\Re(\det(I-\lt_s^n))=\mathrm{exp}\left(-\sum_{k=1}^\infty \frac{1}{k} \Re(\mathrm{Tr}( \lt_s^{nk}))\right) \cos\left( 
\sum_{k=1}^\infty \frac{1}{k} \Im(\mathrm{Tr}( \lt_s^{nk}) \right),$$
and the trace formula combined with the pressure formula (\ref{pressure1}) shows that for all $\epsilon>0$ and $p$ large,
$$\vert \mathrm{Tr}(\lt_s^p)\vert \leq C e^{p(P(\Re(s))+\epsilon)}.$$
We therefore obtain 
$$\vert \Re(Z_\Gamma^{(n)}(2+it))\vert\geq \mathrm{exp}\left({-C}e^{n(P(2)+\epsilon)} \right)\vert  \cos\left(Ce^{n(P(2)+\epsilon)}\right)\vert.$$
Since $P(2)<0$ by Bowen's formula, we have a uniform bound from below as long as $n$ is taken large. $\square$

\bigskip
\noindent
To finish the proof of Proposition \ref{Little}, we use the pointwise estimate of Proposition \ref{est1} which combined with Titchmarsh Lemma gives
$$ \sup_{\sigma_0\leq \sigma \leq 1}\vert \mathrm{Arg}(Z_\Gamma^{(n(T))}(\sigma+iT))\vert =O\left( T^{\delta+\nu P(\sigma_0)} \right)
=O(T),$$
as long as we choose $\nu\leq \frac{1-\delta}{P(\sigma_0)}$.  The proof is done, provided that there are no zeros on
$\{\Im(s)=T\}$ and $\{ \Im(s)=T/2 \}$. If $Z_\Gamma^{(n(T))}(s)$ vanishes on these lines, we simply
replace $T$ by some $T<\widetilde{T}\leq T+1$ so that  $Z_\Gamma^{(n(T))}(s)$ does not vanish on $\{\Im(s)=\widetilde{T} \}$ and
$T/2$ by some $T'$ with $T/2-1\leq T'<T/2$ so that  $Z_\Gamma^{(n(T))}(s)$ does not vanish on $\{\Im(s)=T' \}$.  
We then write 
$$M(\sigma,T)\leq \#\{ s \in \mathcal{R}_X\ :\ \Re(s)\geq \sigma\ \mathrm{and}\ T'\leq \Im(s) \leq \widetilde{T} \}$$
$$\leq C\int_{T'}^{\widetilde{T}}\log\vert Z_\Gamma^{(n(T))}(\sigma_0+it)\vert dt+O(T)$$
$$=C\int_{T/2}^T\log\vert Z_\Gamma^{(n(T))}(\sigma_0+it)\vert dt+O(T),$$
by the pointwise estimate of proposition \ref{est1}. 
$\square$

\bigskip \noindent
We now state the main estimate which is the core argument of the paper. Let $\varphi_0 \in C_0^\infty(\R)$ be a smooth compactly supported function
such that $\mathrm{Supp}(\varphi_0)\subset [-1,+1]$, $\varphi_0>0 $ on $(-1,+1)$ and $\int \varphi_0(x)dx=1$. We define a probability measure $\mu_T$ on $\R$ by the formula
$$\int_\R f d\mu_T:=\frac{1}{T}\int_{-\infty}^{+\infty}\varphi_0 \left(\frac{t-2T}{T} \right)f(t)dt. $$
\begin{propo}
 \label{central}
 There exist $\nu>0$, $0<\rho<1$ such that for all $\sigma>\halfdelta$ one can find $\varepsilon(\sigma)>0$ such that
 $$\sum_{\ell=1}^{N(h)}\int_\R \Vert \lt_{\sigma+it}^{(n(T))}(\mathbf{e}_k^\ell) \Vert^2_{(h)} d\mu_T(t)\leq C \rho^kT^{\delta-\varepsilon(\sigma)},$$
 where we have taken $n(T):=[\nu\log T]$, $h=T^{-1}$.
\end{propo}
\noindent 
The proof of this proposition is postponed to $\S 6$ and occupies the rest of the paper. Let us show how the combination of Proposition \ref{central} and
Proposition \ref{Little} implies Theorem \ref{mainresult}. By Proposition \ref{Little}, we have for  $\halfdelta< \sigma_0<\sigma$,
$$M\left(\sigma,\frac{5}{2}T \right)\leq C\int_{\frac{5}{4}T}^{\frac{5}{2}T}\log\vert Z_\Gamma^{(n(T))}(\sigma_0+it)\vert dt+O(T),$$
and we write
$$\int_{\frac{5}{4}T}^{\frac{5}{2}T} \log\vert Z_\Gamma^{(n(T))}(\sigma_0+it)\vert dt$$
$$\leq \left(\inf_{[-3/4,1/2]}\varphi_0 \right)^{-1} T \int_\R \log\vert Z_\Gamma^{(n(T))}(\sigma_0+it)\vert d\mu_T(t).$$
Using Lemma \ref{det}, we obtain
$$ \int_\R \log\vert Z_\Gamma^{(n(T))}(\sigma_0+it)\vert d\mu_T(t)\leq \sum_{k\in \N}\sum_{\ell=1}^{N(h)}
\int_\R \Vert \lt_{\sigma_0+it}^{(n(T))}(\mathbf{e}_k^\ell) \Vert_{(h)} d\mu_T(t),$$
which by Schwartz inequality is less than
$$Ch^{-\halfdelta} \sum_{k\in \N} \left (    \sum_{\ell=1}^{N(h)}
\int_\R \Vert \lt_{\sigma_0+it}^{(n(T))}(\mathbf{e}_k^\ell) \Vert_{(h)}^2 d\mu_T(t)\right)^{1/2}
\leq C T^{\delta-\epsilon(\sigma_0)/2}.$$
We have therefore obtained
$$M(\sigma,T)=O\left (T^{1+\max\left \{\delta-\half \epsilon(\sigma_0),0 \right \}} \right).$$
To get an estimate on the counting function $N(\sigma,T)$, we just write
$$N(\sigma,T)\leq \sum_{k=0}^{N_0} M\left(\sigma,\frac{T}{2^k}\right)+O(1)$$
where $N_0$ is such that $\frac{T}{2^{N_0}}\leq 1$.  Since we have for $T$ large,
$$\sum_{k=0}^{N_0} M\left(\sigma,\frac{T}{2^k}\right)\leq C T^{1+\max \left \{\delta-\half \epsilon(\sigma_0),0 \right \}},$$
with $C$ independent on $T$, the proof is done.

 \end{section}

\begin{section}{A key lower bound}
Given $\alpha,\beta \in \mathscr{W}_n^j$, we set for all $x\in I_j=\D_j\cap \R$,
$$\Phi_{\alpha,\beta}(x):=\log \vert \gamma'_\alpha(x) \vert -\log \vert \gamma'_\beta(x)\vert .$$
The goal of this section is to prove the following fact which will be a key result in the next section on mean square estimates. 
\begin{propo}
\label{phase}
There exists $\eta_0>0$ and $1>\overline{\theta}>0$ such that for all $n\geq 1$ and all $j=1,\ldots,2p$, we have for all $\alpha\neq \beta \in \mathscr{W}_n^j$,
$$\inf_{x\in I_j} \vert \Phi_{\alpha,\beta}'(x)\vert \geq \eta_0 \overline{\theta}^n.$$
\end{propo}
The proof of this proposition will follow from the next Lemma which is of geometric nature.
\begin{lem}
 Let $\Gamma$ be a Schottky group as above, then there exists a constant $C_\Gamma>0$ such that for all
 $$\gamma\simeq \left ( 
 \begin{array}{cc}
 a&b\\
 c&d
 \end{array}
 \right)\in \Gamma\setminus \{ Id\},\ \vert c \vert \geq C_\Gamma.$$
\end{lem}
\noindent
{\it Proof}. Let $\Gamma\neq Id$ be an element of $\Gamma$, and assume that
$$ \gamma\simeq \left ( 
 \begin{array}{cc}
 a&b\\
 c&d
 \end{array}
 \right),$$
 with $c=0$. By looking at $\gamma^{-1}$ one can assume that $\vert a \vert \geq 1$. In addition, since $\Gamma$ has no non-trivial parabolic element,
 we have actually $\vert a \vert >1$. The action of $\gamma$ on the Riemann sphere is therefore of the type
 $$\gamma(z)=a^2 z+ab.$$
 Now pick an element of $\Lambda(\Gamma)$ which is different from the fixed point of $\gamma$. Its orbit under the action of $(\gamma^n)_{n\in \N}$ goes to infinity as
 $n\rightarrow +\infty$ which contradicts the fact that $\Lambda$ is $\gamma$-invariant and compact. Therefore $c\neq 0$. Remark that by definition of $\Gamma$, we have
 $$\gamma^{-1}(\infty)=-d/c\in \cup_{j=1}^{2p}\mathcal{D}_j.$$ 
 Consequently, one can find a constant $Q_1$ such that  for all $\gamma\neq Id$, $\vert d/c \vert\leq Q_1$. Notice also that we have ($i=\sqrt{-1}$)
 $$\Im(\gamma(i))=\frac{1}{c^2+d^2}.$$
 Since the limit set $\Lambda(\Gamma)$ is compact, the orbit $\Gamma.i$ has to be bounded in $\C$ for the usual euclidean distance. Hence there exists $Q_2$ such that
 for all $\gamma\in \Gamma$,
 $$ \frac{1}{c^2+d^2}\leq Q_2.$$
 As a consequence we get
 $$1\leq c^2(Q_2+Q_1^2),$$
 and the proof is done. $\square$
 
 \bigskip
 \noindent
 We can now complete the proof of Proposition \ref{phase}. Writing
 $$\gamma_\alpha(x)=\frac{a_\alpha x+b_\alpha}{c_\alpha x+d_\alpha},\ \mathrm{with}\ a_\alpha d_\alpha-b_\alpha c_\alpha=1,$$
 we have 
 $$\vert \Phi'_{\alpha,\beta}(x)\vert =2\frac{\vert c_\beta d_\alpha-c_\alpha d_\beta \vert}{\vert c_\beta x+d_\beta\vert \vert c_\alpha x+d_\alpha \vert}
 =2\vert c_\beta d_\alpha-c_\alpha d_\beta \vert (\gamma_\alpha'(x))^{1/2} (\gamma_\beta'(x))^{1/2}.$$
 We can now remark that since $\Gamma$ is a free group, $\alpha\neq \beta$ implies $\gamma_\alpha\circ \gamma_\beta^{-1}\neq Id$, and we have the formula
 $$\gamma_\alpha\circ \gamma_\beta^{-1} \simeq  \left ( 
 \begin{array}{cc}
 a_\alpha&b_\alpha\\
 c_\alpha&d_\alpha
 \end{array}
 \right)
 \left ( 
 \begin{array}{cc}
 d_\beta&-b_\beta\\
 -c_\beta&a_\beta
 \end{array}
 \right)=\left ( 
 \begin{array}{cc}
 *& *\\
 c_\alpha d_\beta-d_\alpha c_\beta&*
 \end{array}
 \right).
$$
We can therefore apply the above Lemma and the proof ends by using bounded distortion and the lower bound for the derivatives.
 \end{section}
 \begin{section}{Mean square estimates}
 The goal of this final section is to prove Proposition \ref{central}. We recall that in the sequel we take 
 $$h=T^{-1},\ n(T)=[\nu \log T],\  0<\nu<<1.$$
 We first start by writing
 $$\sum_{\ell=1}^{N(h)}\int_\R \Vert \lt_{\sigma+it}^{(n(T))}(\mathbf{e}_k^\ell) \Vert^2_{(h)} d\mu_T(t)=
\sum_{j=1}^{2p} \int_{\Omega_j(h)}\int_\R \sum_{\ell=1}^{N(h)}\vert \lt_{\sigma+it}^{(n(T))}(\mathbf{e}_k^\ell)  \vert^2 d\mu_T dm.$$
For all $z\in \D_j$ and $\alpha,\beta \in \mathscr{W}_n^j$ we use the notation
$$\Phi_{\alpha,\beta}(z):=\mathbb{L}(\gamma'_\alpha(z))-\overline{\mathbb{L}(\gamma'_\beta(z))},$$
where $\mathbb{L}$ is the adequate complex logarithm. This notation coincides on the real axis with the one introduced in the previous section.
Remark that since $\Phi_{\alpha,\beta}$ is real on the real axis and because of bounded distortion, we have on any disc $\D_\ell(h)$,
$$\sup_{z\in \D_\ell(h)} \vert e^{iT \Phi_{\alpha,\beta}(z)} \vert=O(1),$$
uniformly on $T$.
For all $z\in \Omega_j(h)$, we have by a change of variable
$$ \int_\R \sum_{\ell=1}^{N(h)}\vert \lt_{\sigma+it}^{(n(T))}(\mathbf{e}_k^\ell)(z)  \vert^2 d\mu_T=$$
$$\sum_{\alpha,\beta \in \mathscr{W}_n^j} (\gamma'_\alpha(z))^\sigma (\overline{\gamma'_\beta(z)})^\sigma G^{(k)}_{\alpha,\beta}(z)
\widehat{\varphi_0}(-T \Phi_{\alpha,\beta}(z)) e^{2iT \Phi_{\alpha,\beta}(z)},$$
where $\widehat{\varphi_0}$ is the usual Fourier transform of $\varphi_0$ defined by
$$\widehat{\varphi_0}(\xi):=\int_{-\infty}^{+\infty} \varphi_0(x) e^{-ix\xi} dx,$$

and $G_{\alpha,\beta}^{(k)}$ is the following sum (see also Lemma \ref{baseestimate}):
$$G_{\alpha,\beta}^{(k)}(z):=\sum_{\ell=1}^{N(h)}  \mathbf{e}_k^\ell \circ \gamma_\alpha(z)   \overline{ \mathbf{e}_k^\ell \circ \gamma_\beta (z)}. $$
The goal now is to gain some decay as $T$ goes to $+\infty$ by using the key observation of Lemma \ref{phase}. We split the above sum into two contributions,
the "diagonal" one plus the "off-diagonal"  one:
$$ \sum_{\alpha,\beta \in \mathscr{W}_n^j} (\gamma'_\alpha(z))^\sigma (\overline{\gamma'_\beta(z)})^\sigma G^{(k)}_{\alpha,\beta}(z)
\widehat{\varphi_0}(-T \Phi_{\alpha,\beta}(z)) e^{2iT \Phi_{\alpha,\beta}(z)}=\mathscr{S}_{diag}+\mathscr{S}_{offdiag},$$
where we have set 
$$\mathscr{S}_{diag}:= \sum_{\alpha \in \mathscr{W}_n^j} \vert \gamma'_\alpha(z) \vert^{2\sigma} G^{(k)}_{\alpha,\alpha}(z)
\widehat{\varphi_0}(-T \Phi_{\alpha,\alpha}(z)) e^{i2T\Phi_{\alpha,\alpha}(z)},$$
$$\mathscr{S}_{offdiag}:=\sum_{\alpha \neq \beta \in \mathscr{W}_n^j} (\gamma'_\alpha(z))^\sigma (\overline{\gamma'_\beta(z)})^\sigma G^{(k)}_{\alpha,\beta}(z)
\widehat{\varphi_0}(-T \Phi_{\alpha,\beta}(z)) e^{2iT \Phi_{\alpha,\beta}(z)}.$$
We first deal with the diagonal contribution, which by bounded distortion and the pressure estimate (together with Lemma \ref{baseestimate}) gives
$$\vert \mathscr{S}_{diag}\vert \leq C e^{n(T)P(2\sigma)}h^{-2} \rho^k, $$
and therefore
$$ \int_{\Omega_j(h)}\vert \mathscr{S}_{diag}\vert dm \leq C h^{-\delta} e^{n(T)P(2\sigma)}\rho^k.$$ 
We clearly have a gain as long as $P(2\sigma)<0$, which by Bowen's formula \cite{Bowen1} amounts to say that $\sigma>\halfdelta$.
To deal with the off-diagonal sum, we will of course use the estimate for $\xi \in \C$ and all $q\geq 0$,
$$ \vert \widehat{\varphi_0}(\xi) \vert \leq C_q \frac{e^{\vert \Im(\xi) \vert}}{(1+\vert \xi \vert)^q},$$
which implies that for $\alpha\neq \beta$ and all $q$, we have
$$\vert \widehat{\varphi_0}(-T \Phi_{\alpha,\beta}(z)) \vert =O_q\left(  (1+\vert T\Phi_{\alpha,\beta}(z)  \vert)^{-q} \right).$$
The trouble comes from $\Phi_{\alpha,\beta}(z)$ which may be vanishing on some part of $\Omega_j(h)$. We prove the following Lemma.
\begin{lem}
 Fix an $\epsilon>0$. There exists $\nu$ small enough (recall that $n(T)=[\nu \log T]$) such that for all $\alpha,\beta \in \mathscr{W}_n^j$ with $\alpha\neq \beta$, one can split the set
 $\mathscr{E}_j(h)=\mathscr{E}_j' \sqcup \mathscr{E}''_j$ so that:
 \begin{itemize}
 \item We have $\# \mathscr{E}_j'=O(h^{-\halfdelta})$.
 \item For all $\ell \in \mathscr{E}_j''$ and $z \in \D_\ell(h)$,
 $$\vert  \Phi_{\alpha, \beta}(z) \vert \geq C h^{1-\halfdelta+\epsilon}. $$
 \end{itemize}
\end{lem}
\noindent
{\it Proof}. By Proposition \ref{phase}, for all $\alpha \neq \beta$, the function $x\mapsto \Phi_{\alpha,\beta}(x)$ is strictly monotonic on the interval $I_j$, and its derivative
is uniformly bounded from below by $C_1 (\overline{\theta})^n$. Two cases can occur: 
\begin{itemize}
\item Either $\Phi_{\alpha,\beta}(x)$ vanishes at some point $x_0\in I_j$ and we set $J_{\alpha,\beta}=[x_0-h^\eta,x_0+h^\eta]$, where $0<\eta<1$ will be adjusted
later on. Then for all
$x\in I_j\setminus J_{\alpha,\beta}$ we have $$\vert \Phi_{\alpha,\beta}(x) \vert\geq C_1 h^\eta (\overline{\theta})^n.$$
\item The map $x\mapsto \Phi_{\alpha,\beta}(x)$ does not vanish on $I_j:=(a_j,b_j)$. According to the sign of $\Phi_{\alpha,\beta}$, we set
$J_{\alpha,\beta}=[a_j,a_j+h^\eta]$ or $J_{\alpha,\beta}=[b_j-h^\eta,b_j]$. In both cases we have for $x\in I_j\setminus J_{\alpha,\beta}$, 
$$\vert \Phi_{\alpha,\beta}(x) \vert\geq C_1 h^\eta (\overline{\theta})^n.$$
\end{itemize}
Now define $\mathscr{E}_j''$ by 
$$\mathscr{E}_j''=\{ \ell \in \mathscr{E}_j(h)\ :\ \D_\ell(h)\cap J_{\alpha,\beta}=\emptyset \}.$$
Given $\ell \in \mathscr{E}_j''$, we have by bounded distortion for all $z\in \D_\ell(h)$,
$$\vert \Phi_{\alpha,\beta}(z) \vert \geq C_1h^\eta (\overline{\theta})^n-C_2 h\geq C_3 h^{\eta+\epsilon},$$
provided $\eta+\epsilon<1$ and $\nu \vert \log \overline{\theta} \vert \leq \epsilon$. It remains to count
$$\mathscr{E}_j':=\mathscr{E}_j(h)\setminus \mathscr{E}''_j.$$
Except for possibly two of them, $\ell \in \mathscr{E}_j'$ implies $\R\cap\D_\ell(h)\subset J_{\alpha,\beta}$, hence by volume comparison
$$(\#\mathscr{E}_j'-2)Ch\leq \vert J_{\alpha,\beta}\vert=h^\eta.$$
Therefore $\# \mathscr{E}_j'=O(h^{\eta-1})$. The proof ends by choosing $\eta=1-\halfdelta$. $\square$

\bigskip \noindent
Going back to $\mathscr{S}_{offdiag}$, we have using Lemma \ref{baseestimate},
$$ \int_{\Omega_j(h)}\vert \mathscr{S}_{offdiag} \vert dm\leq C_q\rho^k h^{-2}\sum_{\alpha\neq \beta} \Vert \gamma'_\alpha \Vert_\infty^\sigma 
\Vert \gamma'_\beta \Vert_\infty^\sigma \int_{\Omega_j(h)} \frac{dm(z)}{(1+\vert T \Phi_{\alpha,\beta}(z) \vert)^q}.$$
We then write
$$\int_{\Omega_j(h)} \frac{dm(z)}{(1+\vert T \Phi_{\alpha,\beta}(z) \vert)^q}=\int_{\Omega_j'(h)} \frac{dm(z)}{(1+\vert T \Phi_{\alpha,\beta}(z) \vert)^q} 
+\int_{\Omega_j''(h)} \frac{dm(z)}{(1+\vert T \Phi_{\alpha,\beta}(z) \vert)^q},$$
with the notations 
$$\Omega_j'(h):=\bigcup_{\ell \in \mathscr{E}_j'} \D_\ell (h)\ \mathrm{and}\  \Omega_j''(h):=\bigcup_{\ell \in \mathscr{E}_j''} \D_\ell (h).$$
We recall to the reader that
$$\# \mathscr{E}_j'=O(h^{-\halfdelta})\  \mathrm{and}\ \# \mathscr{E}_j''=O(h^{-\delta}).$$
Applying the above Lemma and taking $\epsilon>0$ small enough so that $$\halfdelta-\epsilon>0,$$ we get with $T=h^{-1}$ and $q$ large
$$\int_{\Omega_j(h)} \frac{dm(z)}{(1+\vert T \Phi_{\alpha,\beta}(z) \vert )^q}=O\left( h^{2-\halfdelta}\right)+O\left(h^{2-\delta+q(\halfdelta-\epsilon)} \right)$$
$$=O\left(h^{2-\halfdelta} \right).$$
Therefore,
$$\int_{\Omega_j(h)}\vert \mathscr{S}_{offdiag} \vert dm\leq C\rho^k e^{2nP(\sigma)}h^{-\halfdelta}.$$
Adding all of our estimates, we have reached
$$\sum_{\ell=1}^{N(h)}\int_\R \Vert \lt_{\sigma+it}^{(n(T))}(\mathbf{e}_k^\ell) \Vert^2_{(h)} d\mu_T(t)
\leq C\rho^k\left( T^{\delta+\nu P(2\sigma)}+T^{\halfdelta+2\nu P(\sigma)}   \right).$$
Now taking 
$$\nu\leq \min \left \{\frac{\delta}{8\vert \log \overline{\theta} \vert}, \frac{\delta}{16P(\halfdelta)},\frac{1-\delta}{P(\halfdelta)} \right \}$$
yields for all $\sigma \geq \halfdelta$,
$$\sum_{\ell=1}^{N(h)}\int_\R \Vert \lt_{\sigma+it}^{(n(T))}(\mathbf{e}_k^\ell) \Vert^2_{(h)} d\mu_T(t)
\leq C\rho^k\left( T^{\delta+\nu P(2\sigma)}+T^{3\delta/4}   \right).$$
The proof of Proposition \ref{central} is done. 
We can now say a few more words on the function $\tau(\sigma)$ of the main Theorem \ref{mainresult}: using the above choice of $\nu$
then for all $\sigma>\sigma_0\geq \halfdelta$, we can take
$$\tau(\sigma)=\max \left \{ \delta+\frac{\nu}{2} P(2\sigma_0), \frac{7\delta}{8} \right \},$$
so for example taking $\sigma_0=\half (\sigma+\halfdelta)$ gives
\begin{equation}
\label{taufunction}
\tau(\sigma)=\max \left \{ \delta+\frac{\nu}{2} P(\sigma+\halfdelta), \frac{7\delta}{8} \right \}.
\end{equation}
We have not attempted to optimize the choice of constants at all. 
Formula (\ref{taufunction}) shows in addition that $\tau(\sigma)$ is strictly decreasing and convex in a right neighborhood of $\halfdelta$ since the topological
pressure has this property. Moreover, the pressure functional $\sigma \mapsto P(\sigma)$ is real analytic and its derivative at $\sigma=\delta$ is given by
$$P'(\delta)=-\int_{I} \log\vert T' \vert d\mu_\delta<0,$$
where $I=\cup_{i=1}^{2p} I_i$, $T$ is the Bowen-Series map and $\mu_\delta$ is the equilibrium state at $\delta$. In particular, this gives a formula
for  the derivative 
$$\tau'(\halfdelta)=-\frac{\nu}{2}\int_{I} \log\vert T' \vert d\mu_\delta$$
which is of course non-vanishing.

 \end{section}

\end{document}